\documentclass[a4paper,11pt]{article}
\usepackage[english]{babel}
\usepackage{amsmath,amsfonts,latexsym}
\usepackage{amsthm}
\usepackage{amssymb}
\usepackage[dvips]{graphicx}
\usepackage{array}

\makeatletter

\newcommand{\C}{\mathbb C}
\newcommand{\R}{\mathbb R}

\newcommand{\N}{\mathbb N}

\newcommand{\pr}{{\noindent \textsc{Proof}. \, }} 
\newcommand{\epr}{\hspace*{0.1cm}\hfill $\square$\par}

\theoremstyle{plain}
\newtheorem{Thm}{Theorem}
\newtheorem{Def}[Thm]{Definition}
\newtheorem{Pro}[Thm]{Proposition}
\newtheorem{Cor}{Corollary}[Thm]

\theoremstyle{remark}
\newtheorem{Ex}[Thm]{\sc Example}
\newtheorem{Rem}[Thm]{\sc Remark}

\newcommand{\ps@taiplain}{%
     \renewcommand{\@oddhead}{\centerline{\sf {\small The generalized gamma functions}}}%
     \renewcommand{\@evenhead}{\@oddhead}%
     \renewcommand{\@oddfoot}{\centerline{\thepage}}
     \renewcommand{\@evenfoot}{\@oddfoot}}
\makeatother     

\begin{document}

\thispagestyle{plain}
\title{\bf THE GENERALIZED GAMMA FUNCTIONS}

\vspace{-2cm}
\author{ {\sf Tran Gia Loc}
\thanks{The Teacher Training College of Dalat, 29 Yersin, Dalat, Vietnam / {\sc email:} gialoc@gmail.com} \ \   and \ 
{\sf Trinh Duc Tai 
\thanks{The Teacher Training College of Dalat, 29 Yersin, Dalat, Vietnam /  {\sc email:} tductai@gmail.com} }}
\vspace{-2cm}

\date{ }

\maketitle
\pagestyle{taiplain}
\begin{abstract} 
In this paper, we introduce a way to generalize the Euler's gamma function as well as some related special functions. With a given polynomial in one variable $f(t)\ge 0$ , we can associate  a function, so-called "gamma function associated  with $f$", defined by $\Gamma_f(s) := \int_0^\infty f^{s-1} e^{-t} dt$. This function has  many features similar to the Euler's gamma function. We also present some initial results on the  gamma-type functional equation for $\Gamma_f(s)$ in some special cases.
\end{abstract}
\vskip 0.1 cm
{\bf Keywords and phrases:} Gamma function, Bernstein-Sato polynomial.
\vskip 0.2 cm

\section{Introduction}
The Euler's gamma function is one of the most important special functions, because of its role in various fields. The generalization of this famous function has attracted much attention from many mathematicians and physicists. There are some remarkable achievements. 

Barnes has introduced the  multiple gamma function by generalizing the representation formula for Hurwitz's zeta function \cite{Bar, Bar2}.  Post  has given another direction to generalize the gamma function via its limit representation \cite{Pos}. The Vign\'eras' multiple gamma function has been found by virtue of the Bohr-Morellup theorem \cite{Vig}. D\'iaz and Pariguan \cite{Dia} in 2007 introduced the notion of $k-$gamma function by generalizing Pochhammer's symbol. Recently M.Mansour \cite{Man} showed that the  $k-$gamma function can be characterized as the unique solution of a system of functional equations.

Our approach is slightly different from those of the above authors, which is motivated by a question on the existence of  a  functional equation verified by the gamma function associated with a polynomial. For a given polynomial $f(t)$, which is always assumed to be positive when $t > 0$, we define
\begin{equation}\label{E::gammaf}
\Gamma_f(s) := \int_0^\infty f(t)^{s-1}e^{-t}dt  \ .
\end{equation}
The right-hand side of \eqref{E::gammaf} is a holomorphic function in the half of complex plane
 $\Re(s) > 1- \frac{1}{k}$, where $k$ is the multiplicity of $f$ at  $t=0$.
When $f(t) =t$, $\Gamma_f(s)$ is nothing but the Euler's gamma function $\Gamma(s)$. It is well-known that 
$\Gamma(s)$ satisfies the following functional equation
\begin{equation}\label{pth}
\Gamma(s+1) = s\Gamma(s)\ .
\end{equation}
A natural question\footnote{Which is posed to us by Prof. Le Dung Trang.} is that if the gamma function associated $\Gamma_f(s)$ verifies some kind of   functional equation, type of \eqref{pth}. More concretely, it is expected  that there exists a polynomial $\mathcal{B}(s)$ such that the following is true
\begin{equation}\label{E::ptham}
\Gamma_f(s+1) = \mathcal{B}(s)\Gamma_f(s)\ .
\end{equation}
In general, $\mathcal{B}(s)$ should be depend on $f$ and is conjectured to be the Bernstein-Sato polynomial of $f$.
In section 2, we give a positive answer for this question in a very special case where $f$ is a monomial. There we also present  some interesting properties of $\Gamma_f(s)$, in comparison with those in \cite{Dia}. Section 3 is devoted to define zeta and  beta functions associated with $f$.  We give a counter-example for the conjecture in the last section.

\section{The gamma function associated with $f(t) = t^k$}
As usual,  a power of complex variable $t^z$ is defined by $t^z := e^{z\ln t}$, where $\ln t$ is the principal value of  the logarithm.

\begin{Def}
$$\Gamma_{t^k}(s) := \int_0^\infty t^{k(s-1)} e^{-t} dt \quad , \quad Res > 1-\frac{1}{k} \ .$$   
\end{Def}
From the above definition, it follows that
\begin{Cor}\ 
\begin{description}
\item[\qquad (i)] $\Gamma_{t^k}(1) = 1$,
\item[\qquad (ii)] $\Gamma_t(s) = \Gamma(s)$,
\item[\qquad (iii)] $\Gamma_{t^k}(s) = \Gamma[k(s-1)+1] $.
\end{description}
\end{Cor}

\subsection*{Bernstein-Sato polynomials} 
Let K be a field of characteristic zero and $s$ be a parameter. We recall here the so-called the Bernstein-Sato polynomial (see \cite{Gra}, p. 235)
\begin{Thm}
Let $f\in K[x_1,\cdots,x_n]$ be a non zero polynomial. There is a polynomial $\mathcal{B}(s)\in K[s]$ and a differential operator $P(x,\dfrac{\partial}{\partial x},s) \in A_n(K)[s]$ such that 
\begin{equation}\label{E::ptB}
P(x,\dfrac{\partial}{\partial x},s)f^s = \mathcal{B}(s)f^{s-1} \ . 
\end{equation}
\end{Thm} 

The set of  polynomials $\mathcal{B}(s)$ veryfying the equation \eqref{E::ptB}  is clearly an ideal of $K[s]$, 
which is a principal ideal. The Bernstein-Sato polynomial of $f$ is by definition the monic generator of this ideal. Let us denote it  by $b_f(s)$ or simply $b(s)$.

\begin{Ex} We consider a simple example with $f(t) = t^k, \  k\in\N,\ k>0$. Then
$$\dfrac{d}{dt}f^s = ks t^{k-1}f^{s-1}\ . $$
Assume that $ \dfrac{d^m}{dt^m}f^s = C_m(s)t^{k-m}f^{s-1}.$ Hence
$$\begin{array}{lllll}
\dfrac{d^{m+1}}{dt^{m+1}}f^s   & = & C_m(s)(k-m)t^{k-m-1}f^{s-1} + C_m(s)(s-1)kt^{2k-m-1}f^{s-2} \\
\                                                          & = & C_m(s)\Big((k-m) + (s-1)k\Big)t^{k-m-1}f^{s-1} \\
\end{array}$$
Therefore
\begin{equation}\label{E::ptdq}
C_{m+1}(s) = C_m(s)(ks-m) \ .
\end{equation}      
The formula \eqref{E::ptdq} gives
$$ \begin{array}{llll}
C_1(s) & = & ks \\
C_2(s) & = & ks(ks-1) \\
\        & \cdots & \\
\end{array}$$
$$  C_k(s) = ks(ks-1)(ks-2)\cdots (ks-(k-1))\ .$$  
Therefore, 
$$P\Big(t,\dfrac{d}{dt},s\Big)\big(t^k\big)^s = \mathcal{B}(s)\big(t^k\big)^{s-1} \ ,$$
where $P\Big(t,\dfrac{d}{dt},s\Big) = \dfrac{d^k}{dt^k}$ and 
$$ \mathcal{B}(s) := ks(ks-1)(ks-2)\cdots (ks-(k-1))\ .$$
As a consequence,
the corresponding Bernstein-Sato polynomial is $b(s) = s(s-1/k)(s-2/k)\cdots (s-(k-1)/k)$
\end{Ex}

\subsection{The properties of  $\Gamma_{t^k}$}
\begin{Pro}\label{ttt}
\begin{equation}\label{E::ptham1}
\Gamma_{t^k}(s+1) = \mathcal{B}(s)\Gamma_{t^k}(s) \ .
\end{equation}
\end{Pro}
\pr
Let $\mathcal{L}$ be Laplace transform which is defined by  
$$\mathcal{L}\{f(t);\alpha\} := \int_0^\infty f(t) e^{-\alpha t}dt \quad,\quad Re\alpha >0 \ .$$
Then  
$$\Gamma_{t^k}(s) = \mathcal{L}\big\{t^{k(s-1)};1\big\} \ .$$
By \cite{Deb}, p.144 we have  
\begin{equation}\label{E::qhe}
\mathcal{L}\big\{f^{(k)}(t);\alpha\big\} = \alpha^n\mathcal{L}\{f(t);\alpha\} - \alpha^{n-1}f(0) -\alpha^{n-2}f'(0) - \cdots - \alpha f^{(n-2)}(0) - f^{(n-1)}(0) \ .
\end{equation}
For $f(t) = t^{ks}$, we have 
\begin{equation}\label{E::ptham2}
f^{(k)}(t) = \dfrac{d^k}{dt^k}(t^k)^s = \mathcal{B}(s)(t^k)^{s-1} \ .
\end{equation}
From \eqref{E::qhe} and \eqref{E::ptham2}, we have  
$$\int_0^\infty t^{ks} e^{-t}dt = \int_0^\infty \mathcal{B}(s)(t^k)^{s-1} e^{-t} dt = \mathcal{B}(s)\int_0^\infty (t^k)^{s-1} e^{-t} dt \ .$$
It follows \eqref{E::ptham1} is true.
\epr
\vskip 0.3cm
\begin{Rem} 
The functional equation \eqref{E::ptham1} can be put in the form 
$$ \Gamma_{t^k}(s)=\dfrac{\Gamma_{t^k}(s+1)}{\mathcal{B}(s)}$$
 By itering this process, we obtain
$$ \Gamma_{t^k}(s) = \dfrac{\Gamma_{t^k}(s+n)}{\mathcal{B}(s)\cdots \mathcal{B}(s+(n-1))} \ .$$
\end{Rem}
\begin{Cor} The $\Gamma_{t^k}(s)$ admits an analytic continuation as a meromorphic function with poles on the set :
$$ \{s+k \ | \ s\in\C, k\in\N,\ \mathcal{B}(s+k) = 0\}.$$
\end{Cor}

\begin{Pro}
$$\Gamma_{t^k}(s) = \lim_{n\to\infty}\dfrac{n! n^{k(s-1)+1}}{[k(s-1)+1][k(s-1)+2]\cdots [k(s-1)+(n+1)]} \ .$$
\end{Pro}
\pr
We have 
$$\Gamma_{t^k}(s) = \int_0^\infty t^{k(s-1)} e^{-t}dt = \lim_{n\to\infty}\int_0^n t^{k(s-1)}\Big(1- \dfrac{t}{n}\Big)^n dt \ .$$
By setting $\tau = \frac{t}{n}$, it follows that
$$\Pi (s,n) := \int_0^n t^{k(s-1)}\Big(1- \dfrac{t}{n}\Big)^n dt = n^{k(s-1)+1}\int_0^1 \tau^{k(s-1)}(1-\tau)^n d\tau \ .$$ 
By integration by parts, we have 
\begin{align*}
\int_0^1  \tau^{k(s-1)}(1-\tau)^n d\tau  &= \dfrac{\tau^{k(s-1)+1}(1-\tau)^n}{k(s-1)+1}\Big |_0^1 + \dfrac{n}{k(s-1)+1}
                                                                             \int_0^1 \tau^{k(s-1)+1}(1-\tau)^{n-1} d\tau \\ 
\\
                                                                   & = \dfrac{n(n-1)}{[k(s-1)+1][k(s-1)+2]}\int_0^1 \tau^{k(s-1)+2}(1-\tau)^{n-2}d\tau \notag \\
                                                                  & \cdots  \\
                                                                  & =   \dfrac{n(n-1)\cdots 2.1}{[k(s-1)+1][k(s-1)+2]\cdots [k(s-1)+n]}\int_0^1 \tau^{k(s-1)+n}d\tau \\
\\
                                                                  &= \dfrac{1.2 \cdots (n-1)n}{[k(s-1)+1][k(s-1)+2]\cdots [k(s-1)+(n+1)]} \ .
\end{align*}
So 
$$\Pi(s,n) = \dfrac{1.2 \cdots (n-1)n}{[k(s-1)+1][k(s-1)+2]\cdots [k(s-1)+(n+1)]}. n^{k(s-1)+1} $$
and
$$\Gamma_{t^k}(s) = \lim_{n\to\infty}\Pi(s,n) \ .$$ 
\epr

\begin{Pro}
\begin{equation}\label{E::gamma1}
\dfrac{1}{\Gamma_{t^k}(s)} = [k(s-1)+1] e^{\gamma[k(s-1)+1]} \prod_{n=1}^\infty \Big(1+\frac{k(s-1)+1}{n}\Big) e^{-\frac{k(s-1)+1}{n}} \ ,
\end{equation}
where $\gamma = \lim_{n\to\infty}\Big(1 +\frac{1}{2} + \cdots + \frac{1}{n} - logn\Big)$ .
\end{Pro}
\pr
We have 
\begin{align*}
\Pi(s,n)  & =  \dfrac{e^{[k(s-1)+1]logn}}{[k(s-1)+1][1 + \frac{k(s-1)+1}{1}][1+\frac{k(s-1)+1}{2}]\cdots [1 + \frac{k(s-1)+1}{n}]} \\
\\
              & =  \dfrac{e^{[k(s-1)+1]\big(logn -1 - \frac{1}{2} - \cdots - \frac{1}{n}\big)} e^{\big(\frac{k(s-1)+1}{1} + \frac{k(s-1)+1}{2} + \cdots + \frac{k(s-1)+1}{n}\big)}}{[k(s-1)+1][1 + \frac{k(s-1)+1}{1}][1+\frac{k(s-1)+1}{2}]\cdots [1 + \frac{k(s-1)+1}{n}]} \ . \\
\intertext{Set $\gamma_n = (1 + \frac{1}{2} + \cdots + \frac{1}{n} - logn)$,} 
              & = \dfrac{e^{-\gamma_n [k(s-1)+1]}}{k(s-1)+1}\prod_{i=1}^n \dfrac{e^{\frac{k(s-1)+1}{i}}}{1+\frac{k(s-1)+1}{i}}
\end{align*}
It follows that
$$\Gamma_{t^k}(s) = \lim_{n\to\infty}\Pi(s,n) =  \dfrac{e^{-\gamma [k(s-1)+1]}}{k(s-1)+1}\prod_{n=1}^\infty \dfrac{e^{\frac{k(s-1)+1}{n}}}{1+\frac{k(s-1)+1}{n}} \ .$$
This completes the proof.
\epr
\vskip 0.5cm
\begin{Pro}
\begin{equation}\label{gamma-sin}
\Gamma_{t^k}(s)\Gamma_{t^k}(1-s) = \dfrac{\pi}{\sin(\pi ks)} \prod_{i=1}^{k-1}\dfrac{1}{k(s-1)+i} \ .
\end{equation}
\end{Pro}
\pr
We have
$$\Gamma_{t^k}(s) = \int_0^\infty t^{k(s-1)} e^{-t}dt = \Gamma[k(s-1) +1] \quad \text{and} \quad \Gamma_{t^k}(1-s) = \Gamma(1-ks) .$$
On the other hand, we have the well-known functional equation $\Gamma(s)\Gamma(1-s) = \dfrac{\pi}{\sin (\pi s)}$, therefore 
$$\Gamma(ks)\Gamma(1-ks) = \dfrac{\pi}{\sin[\pi(ks)]}\ .$$
But
\begin{align*}
\Gamma_{t^k}(s) & = \Gamma\big[k(s-1)+1 \big]  = \dfrac{\Gamma \big[k(s-1)+2\big]}{[k(s-1)+1]} =  \dfrac{\Gamma \big[k(s-1)+3\big]}{[k(s-1)+1][k(s-1)+2]} \\ 
                                           & = \cdots \\
                                           & = \dfrac{\Gamma \big[k(s-1)+k \big]}{[k(s-1)+1][k(s-1)+2]\cdots {[k(s-1)+(k-1)]}}\\ \\
                                           & = \Gamma (ks) \prod_{i=1}^{k-1}\dfrac{1}{[k(s-1)+i]}
\end{align*}                                           
So
$$\Gamma_{t^k}(s) \Gamma_{t^k}(1-s) = \Gamma (ks)\Gamma(1-ks) \prod_{i=1}^{k-1}\dfrac{1}{[k(s-1)+i]} 
                                                                         = \dfrac{\pi}{\sin\big(\pi ks\big)} \prod_{i=1}^{k-1}\dfrac{1}{k(s-1)+i} \ .$$
This completes the proof.
\epr
\begin{Rem} Throught the above results, we can conclude that the gamma function associated to $f(t) = t^k$ 
has almost properties similar to the Euler's gamma function. In particular, the functional equation \eqref{E::ptham} is true in this case.
\end{Rem}

\subsection{Asymptotic expansion of  $\Gamma_{t^k}$}
We recall here a classical result on asymptotic expansion which can be found in \cite{Etin}.
\begin{Thm}
Assume that $f: (a,b) \longrightarrow \R$, with $a,\ b \in [0,+\infty)$ attains a global minimum at a unique point $c\in (a,b)$, such that $f^{''}(c)>0$. Then one has 
\begin{equation}\label{E::tiemcan}
\int_a^b g(x)e^{-\frac{f(x)}{h}}dx = h^{\frac{1}{2}} e^{-\frac{f(c)}{h}}\sqrt{2\pi}\dfrac{g(c)}{\sqrt{f^{''}(c)}} + O(h) \ .
\end{equation}
\end{Thm}
The below proposition gives us the asymptotic behavior  of  $\Gamma_{t^k}$
\begin{Pro}
For $Res>0$, the following identity holds 
$$\Gamma_{t^k}(s) = \dfrac{(2\pi)^{\frac{1}{2}}(ks)^{ks+\frac{1}{2}}}{\mathcal{B}(s)} e^{-ks} + O\Big(\frac{1}{s}\Big) \ .$$
\end{Pro}

\pr
We have  
$$\Gamma_{t^k}(s+1) = \int_0^\infty t^{ks} e^{-t}dt.$$
By making the change of variable $t = s\omega$ we have
$$\Gamma_{t^k}(s+1) = s^{ks+1}\int_0^\infty \omega^{ks} e^{-s\omega}d\omega =  s^{ks+1}\int_0^\infty e^{-s(\omega - k\log\omega)} d\omega \ .$$
Let $f(\omega) = \omega - k\log\omega$. Clearly $f^{'}(\omega) =0$ if and only if $\omega = k$. On the other hand $f^{''}(k) = k^{-1}>0$. 
From \eqref{E::tiemcan}, we have
$$\int_0^\infty e^{-s(\omega - k\log\omega)} d\omega = \Big(\dfrac{1}{s}\Big)^\frac{1}{2}. e^{-\frac{k-k\log k}{\frac{1}{s}}}\sqrt{2\pi}\dfrac{1}{\sqrt{k^{-1}}} 
                         + O\Big(\frac{1}{s}\Big)     = \dfrac{(2\pi)^{\frac{1}{2}}}{s^{\frac{1}{2}}} k^{ks+\frac{1}{2}}e^{-ks}+O\Big(\frac{1}{s}\Big)\ .$$
Therefore
\begin{align*}
\Gamma_{t^k}(s+1)    &= (2\pi)^{\frac{1}{2}} s^{ks+\frac{1}{2}} k^{ks+\frac{1}{2}} e^{-ks} + O\Big(\frac{1}{s}\Big) 
\end{align*}
By virtue of \eqref{E::ptham1}, we complete the proof.
\epr

 \subsection{The relation between $\Gamma_{t^k}$ and $\Gamma_k$ }
 The family of function $\Gamma_k$ ($k > 0$) , which is called $k-$gamma function, is defined by (see \cite{Dia})
$$\Gamma_k(s) := \int_0^\infty t^{s-1}e^{-\frac{t^k}{k}}dt \ ,$$
where $s\in \C, \ Res > 0$. The following proposition shows the closely relation between our gamma funtion $\Gamma_{t^k}$ and this $k-$gamma function.

\begin{Pro}
\begin{equation}\label{E::ptham3}
\Gamma_{t^k}(s) = \dfrac{k^{ks}s}{\mathcal{B}(s)}\Gamma_{\frac{1}{k}}(s) \ .
\end{equation}
\end{Pro}

\pr
We have 
$$\Gamma_{t^k}(s) = \int_0^\infty t^{k(s-1)}e^{-t}dt \ .$$
By making the change of variable $t=k \omega^{\frac{1}{k}}$. Then 
\begin{align*}
\Gamma_{t^k}(s) & = k^{k(s-1)}\int_0^\infty \omega^{\big(s-1+\frac{1}{k}\big)-1} e^{-\big(\frac{\omega^{\frac{1}{k}}}{\frac{1}{k}}\big)} d\omega  \\
                                & = k^{k(s-1)}\Gamma_{\frac{1}{k}}\Big(s-1+\frac{1}{k}\Big)
\end{align*}                                
By replacing s with s+1, we have 
\begin{equation}\label{E::ptham4}
\Gamma_{t^k}(s+1) = k^{ks}\Gamma_{\frac{1}{k}}\Big(s+\frac{1}{k}\Big) \ .
\end{equation}
On the other hand, by \cite{Dia}, p. 183 
\begin{equation}\label{E::ptham5}
\Gamma_k(s+k) = s\Gamma_k(s) \ .
\end{equation}
From \eqref{E::ptham1}, \eqref{E::ptham4} and \eqref{E::ptham5}, it follows that 
$$\mathcal{B}(s)\Gamma_{t^k}(s) = k^{ks}s. \Gamma_{\frac{1}{k}}(s) \ .$$
This completes the proof. 
\epr
\vskip 0.5cm
\begin{Pro}
$$\Gamma_{t^k}(s).\Gamma_{t^k}\Big(\frac{1}{k}-s\Big) = \dfrac{s(1-ks)}{\mathcal{B}(s)\mathcal{B}\big(\frac{1}{k}-s\big)}\ .\ \dfrac{\pi}{\sin(\pi ks)}\ .$$
\end{Pro}

\pr
By \cite{Dia}, p. 183 
$$\Gamma_k(s)\Gamma_k(k-s) = \dfrac{\pi}{\sin\big(\frac{\pi s}{k}\big)} \ .$$
By replacing k with $\frac{1}{k}$, we have
$$\Gamma_{\frac{1}{k}}(s) . \Gamma_{\frac{1}{k}}\Big(\frac{1}{k}-s \Big) = \dfrac{\pi}{\sin(\pi ks)} \ .$$ 
It follows from \eqref{E::ptham3}that
$$\Gamma_{\frac{1}{k}}(s) = \dfrac{\mathcal{B}(s)\Gamma_{t^k}(s)}{k^{ks}s} \qquad \text{and} \qquad  \Gamma_{\frac{1}{k}}\Big(\frac{1}{k} - s\Big) = \dfrac{\mathcal{B}\Big(\frac{1}{k} - s\Big)\Gamma_{t^k}\Big(\frac{1}{k} - s\Big)}{k^{k\big(\frac{1}{k} - s\big)}\Big(\frac{1}{k} - s\Big)} \ .$$
Then
$$\dfrac{\mathcal{B}(s)}{k^{ks}s}.\dfrac{\mathcal{B}\Big(\frac{1}{k} - s\Big)}{k^{1-ks}\Big(\frac{1}{k} - s\Big)}.\Gamma_{t^k}(s).\Gamma_{t^k}\Big(\frac{1}{k} - s\Big) = \dfrac{\pi}{\sin(\pi ks)} \ . $$
Hence, we obtain
$$\dfrac{\mathcal{B}(s)\mathcal{B}\Big(\frac{1}{k} - s\Big)}{s(1-ks)}.\Gamma_{t^k}(s).\Gamma_{t^k}\Big(\frac{1}{k} - s\Big) = \dfrac{\pi}{\sin(\pi ks)} \ . $$
This completes the proof. 
\epr

\section{Generalized Zeta and Beta functions}
\quad In this section we define $f-$Beta and $f-$Zeta functions associated with a polynomial $f$. For $f(t) = t^k$, we prove that 
they have many properties similar  to those of classical Zeta and Beta functions.
\vskip 0.1 cm
The Zeta function is studied first by L. Euler (1707-1783), who considered only real values of $s$. The notion of $\zeta(s)$ as a function of the complex variable $s$ is due to B. Riemann (1826 - 1866). The Riemann zeta function is defined by  
$$\zeta(s) = \sum_{n=1}^\infty \dfrac{1}{n^s}\ ,\qquad Res >1.$$
We have the well-known functional equation 
$$\zeta(s) = \dfrac{1}{\Gamma(s)}\int_0^\infty t^{s-1}(1-e^{-t})^{-1} e^{-t}dt \ .$$ 
Hurwitz's zeta function is defined by 
$$\zeta_H(s,a) := \sum_{n=0}^\infty \dfrac{1}{(n+a)^s} \quad , \quad Res > 1, \ a \neq 0,-1,-2,\cdots $$
This is a generalization of the Riemann zeta function, and we also have the well-known similar functional equation
$$\zeta_H(s,a) = \dfrac{1}{\Gamma(s)}\int_0^\infty t^{s-1}(1-e^{-t})^{-1}e^{-at}dt \ .$$   
\noindent The Beta function or Euler integral of the first kind is defined by 
$$B(p,q) := \int_0^1 t^{p-1}(1-t)^{q-1} dt \quad , \quad Re(p) > 0,\ Re(q) > 0. $$
We have (see \cite{Rai})
$$ B(p,q) = \dfrac{\Gamma(p)\Gamma(q)}{\Gamma(p+q)} \ .$$

\subsection{$f-$Beta and $f-$Zeta functions}
Let $f$ be a polynomial, $f-$Beta and $f-$Zeta functions are defined by
\begin{Def}
$$ B_f(p,q) := \dfrac{\Gamma_f(p)\Gamma_f(q)}{\Gamma_f(p+q)}  \quad , \quad Re(p) > 0,\ Re(q) > 0 \ .$$

$$\zeta_f(s) := \dfrac{1}{\Gamma_f(s)}\int_0^\infty f^{s-1}(1-e^{-t})^{-1} e^{-t} dt\quad, \quad Res>1 \ .$$
\end{Def}
From the above definition, we have the below proposition

\begin{Pro}
\begin{equation}\label{E::f-zeta}
\zeta_f (s) \Gamma_f (s) = \sum_{n=0}^\infty \int_0^\infty f^{s-1} e^{-(n+1)t} dt \ .
\end{equation}
\end{Pro}
\pr
Since $(1-e^{-t})^{-1} = \sum_{n=0}^\infty e^{-nt}$, therefore 
$$\zeta_f(s)\Gamma_f(s) = \int_0^\infty f^{s-1}\Big(\sum_{n=0}^\infty e^{-nt}\Big) e^{-t} dt = \sum_{n=0}^ \infty \int_0^\infty f^{s-1} e^{-(n+1)t}dt . $$
\epr

\subsection{$f-$Beta and $f-$Zeta functions in case $f(t)=t^k$} 
Let $f(t) = t^k$, with $k\in \N, \ k>0$. Then $t^k-$Beta and $t^k-$Zeta are defined by
\begin{equation}\label{E::t^k-beta} 
B_{t^k}(p,q) = \dfrac{\Gamma_{t^k}(p)\Gamma_{t^k}(q)}{\Gamma_{t^k}(p+q)}  \quad , \quad Re(p) > 0,\ Re(q) > 0 \ ,
\end{equation}

$$\zeta_{t^k}(s) := \dfrac{1}{\Gamma_{t^k}(s)}\int_0^\infty t^{k(s-1)}(1-e^{-t})^{-1} e^{-t} dt\quad, \quad Res>1 \ .$$

From \eqref{E::f-zeta} we have the corollary  
\begin{Cor}
\begin{equation}\label{E::t^k-zeta}
\zeta_{t^k}(s)\Gamma_{t^k}(s) = \sum_{n=0}^\infty \int_0^\infty t^{k(s-1)} e^{-(n+1)t} dt \ .
\end{equation}
\end{Cor}

The below proposition gives a relation between function $\zeta_{t^k}$ and Riemann zeta function.

\begin{Pro}
$$\zeta_{t^k}(s) = \zeta\big[k(s-1)+1\big] \ .$$
\end{Pro}
\pr
By making the change of variable $\omega = (n+1)t$ in \eqref{E::t^k-zeta}, we have 
\begin{align*}
\zeta_{t^k}(s)\Gamma_{t^k}(s)  & = \sum_{n=0}^\infty \dfrac{1}{(n+1)^{k(s-1)+1}} \int_0^\infty \omega^{k(s-1)}  e^{-\omega} d\omega \\
                                                          & = \Gamma_{t^k}(s)\sum_{n=0}^\infty \dfrac{1}{(n+1)^{k(s-1)+1}} \\
                                                          & = \Gamma_{t^k}(s)\sum_{n=1}^\infty \dfrac{1}{n^{k(s-1)+1}}\\
                                                          & = \Gamma_{t^k}(s) . \zeta\big[k(s-1)+1\big] \ .
\end{align*}
This proves the proposition.
\epr

\begin{Pro}
$$B_{t^k}(p,q) = \dfrac{kpq}{p+q}.\dfrac{\mathcal{B}(p+q)}{\mathcal{B}(p)\mathcal{B}(q)}.B(kp,kq)\ .$$
\end{Pro}

\pr
From \eqref{E::ptham3} and \eqref{E::t^k-beta}, we have  
$$ B_{t^k}(p,q) = \dfrac{\dfrac{pk^{kp}}{\mathcal{B}(p)}\Gamma_{\frac{1}{k}}(p)\dfrac{qk^{kq}}{\mathcal{B}(q)}
                                  \Gamma_{\frac{1}{k}}(q)}{\dfrac{(p+q)k^{k(p+q)}}{\mathcal{B}(p+q)}\Gamma_{\frac{1}{k}}(p+q)} 
                             = \dfrac{pq}{p+q}.\dfrac{\mathcal{B}(p+q)}{\mathcal{B}(p)\mathcal{B}(q)}.
                                 \dfrac{\Gamma_{\frac{1}{k}}(p)\Gamma_{\frac{1}{k}}(q)}{\Gamma_{\frac{1}{k}}(p+q)}$$
By \cite{Man}, p. 187 
$$ B_k(p,q)=\dfrac{\Gamma_k(p)\Gamma_k(q)}{\Gamma_k(p+q)} \qquad \text{and}\qquad B_k(p,q) 
                    = \dfrac{1}{k}B\Big(\frac{p}{k},\frac{q}{k}\Big) \ .$$
So 
$$B_{t^k}(p,q) =    \dfrac{pq}{p+q}.\dfrac{\mathcal{B}(p+q)}{\mathcal{B}(p)\mathcal{B}(q)}.B_{\frac{1}{k}}(p,q) 
                            =  \dfrac{pq}{p+q}.\dfrac{\mathcal{B}(p+q)}{\mathcal{B}(p)\mathcal{B}(q)}.kB(kp,kq) \ .$$
The proof is complete.
\epr

\section{A functional equation for  $\Gamma_f$ for a quadratic polynomial}
In this section, we consider the quadratic case as a counter example for the truth of \eqref{E::ptham}
Let $f(t) = t^2 + bt +c $. Then  
$$ \begin{array}{llllll}
\dfrac{d}{dt}f^s          & = & s(2t+b)(t^2+bt+c)^{s-1}  \ = \  s(2t+b)f^{s-1} \\
\\
\dfrac{d^2}{dt^2}f^s  & = & 2sf^{s-1}+s(s-1)(2t+b)^2f^{s-2}                 \\
                                        & = & 2sf^{s-1}+s(s-1)(4t^2+4bt+b^2)f^{s-2}\\
                                     & = & 2sf^{s-1}+s(s-1)\Big[4(t^2+4bt+c)+(b^2-4c)\Big]f^{s-2}  \\
                                     & = &  2s(2s-1)f^{s-1} + (b^2-4c)s(s-1)f^{s-2} \\
\end{array}
$$
Therefore
$$ \Big[(t^2+bt+c)\dfrac{d^2}{dt^2} - 2s(2s-1)\Big]f^{s}   =  (b^2-4c)s(s-1)f^{s-1} .$$
So 
$$P\Big(t,s,\dfrac{d}{dt}\Big) = \Big[(t^2+bt+c)\dfrac{d^2}{dt^2} - 2s(2s-1)\Big] \qquad \text{and} \qquad \mathcal{B}(s) =(b^2-4c) s(s-1) .$$
\vskip 0.3cm
\noindent Consider the function
$$ \Gamma_f(s) = \int_0^\infty (t^2+bt+c)^{s-1}e^{-t}dt .$$
We have
\begin{align*}
\mathcal{B}(s)\Gamma_f(s)  & = \int_0^\infty \mathcal{B}(s)(t^2+bt+c)^{s-1}e^{-t}dt \\
								&= \int_0^\infty  \Big[(t^2+bt+c)\dfrac{d^2}{dt^2} - 2s(2s-1)\Big](t^2+bt+c)^{s}e^{-t}dt \\
                                & = \int_0^\infty  (t^2+bt+c)^{s} \Big[(t^2+bt+c)\dfrac{d^2}{dt^2} - 2s(2s-1)\Big]^* e^{-t}dt \qquad\qquad \\
                               & = \int_0^\infty (t^2+bt+c)^s \Big[\dfrac{d^2}{dt^2}(t^2+bt+c)e^{-t} - 2s(2s-1)e^{-t}\Big]dt \qquad \\
                               & =  \int_0^\infty f^s \Big[2e^{-t} -(2t+b)e^{-t} - (2t+b)e^{-t} + (t^2+bt+c)e^{-t} - 2s(2s-1)e^{-t}\Big]dt \\
                               & = \int_0^\infty f^s\Big[(t^2+bt+c) -2(2t+b) + 2-2s(2s-1)\Big]e^{-t} dt \qquad \\
                               & = \int_0^\infty f^{s+1}e^{-t}dt - 2(s-1)(2s+1)\int_0^\infty f^s e^{-t}dt - 2\int_0^\infty (2t+b)(t^2+bt+c)^s e^{-t}dt \\
                              & = \Gamma_f(s+2) - 2(s-1)(2s+1)\Gamma_f(s+1) - 2\int_0^\infty (2t+b)(t^2+bt+c)^s e^{-t} dt .
\end{align*}                             
By integration by parts $\big(u=e^{-t}, \quad dv = (2t+b)(t^2+bt+c)^s dt \big)$, we have   
\begin{align*}
\mathcal{B}(s)\Gamma_f(s)  & = \Gamma_f(s+2) - 2(s-1)(2s+1)\Gamma_f(s+1) - \\
                                &\ \ \qquad - 2\Big[\dfrac{(t^2+bt+c)^{s+1}}{s+1}e^{-t}\Big|_0^\infty + \dfrac{1}{s+1}\int_0^\infty (t^2+bt+c)^{s+1} e^{-t}dt\Big] \\
                                & = \Gamma_f(s+2) - 2(s-1)(2s+1)\Gamma_f(s+1) +\dfrac{2c^{s+1}}{s+1} - \dfrac{2}{s+1}\Gamma_f(s+2).
\end{align*}
So we get the proposition 
\begin{Pro}
\begin{equation}\label{E::pthambachai}
\Big(1-\dfrac{2}{s+1}\Big)\Gamma_f(s+2) = \mathcal{B}(s)\Gamma_f(s) + 2(s-1)(2s+1)\Gamma_f(s+1) + \dfrac{2c^{s+1}}{s+1}.
\end{equation}
\end{Pro}
\vskip 0.2cm
\begin{Rem} The functional equation \eqref{E::pthambachai} is a type of second-order difference equation and impossible to be reduced to that of first order as \eqref{E::ptham}. 
In general, the functional equation \eqref{E::ptham}
is not true for any polynomial. Meanwhile, we guess that for generic polynomials $f(t)$, the gamma function  associated $\Gamma_f(s)$ must satisfy a  difference equation whose order is at most the degree of $f$.
\end{Rem}

\indent {\bf Acknowledgments}. This paper is supported by Vietnam’s National Foundation for Science and Technology Development (NAFOSTED). We would like to thank Professor L\^e D\~ung Tr\'ang for his valuable suggestions on this paper.

\end{document}